\documentclass[12pt]{amsart}
\usepackage{amsmath,amsfonts,amssymb,amsthm}
\usepackage{euscript}
\newcommand{\version}{%
 (February 12, 2006)}
\input xy
\xyoption{all}
\begin{document}
%
\newcommand{\refequ}[1]{$(\ref{#1})$}
\theoremstyle{plain}
    \newtheorem{thm}     {Theorem}[section]
    \newtheorem*{thmA}     {Theorem A}
    \newtheorem*{thmB}     {Theorem B}
    \newtheorem{prop}[thm]    {Proposition}
    \newtheorem{propdef}[thm]    {Proposition-Definition}
    \newtheorem{cor}[thm]     {Corollary}
    \newtheorem{lemma}[thm]   {Lemma}
    \newtheorem{conj}[thm]    {Conjecture}
\theoremstyle{definition}
    \newtheorem{defin}[thm]{Definition}
    \newtheorem{example}[thm]{Example}
    \newtheorem{examples}[thm]{Examples}
    \newtheorem{notation}[thm]{Notation}
\newcommand{\nothm}{\thethm\ \addtocounter{thm}{1}}
\newcommand{\nothmlbl}[1]{\newcounter{#1}\setcounter{#1}{\thethm}\nothm}
\newenvironment{myeq}[1][]
{\stepcounter{thm}\begin{equation}\tag{\thethm}{#1}}
{\end{equation}}
\newcommand{\Hom}{\operatorname{Hom}}
\newcommand{\eHom}{\EuScript Hom}
\newcommand{\parf}{\operatorname{parf}}
\newcommand{\Spec}{\operatorname{Spec}}
\newcommand{\qc}{\operatorname{qc}}
\newcommand{\Ass}{\operatorname{Ass}}
\newcommand{\ann}{\operatorname{ann}}
%
\newcommand{\rmk}{{\bf Remark }}
\newcommand{\rmks}{{\bf Remarks }}
\newcommand{\finrmk}{\medbreak\noindent}
\newcommand{\xcite}[1]{\cite{{#1}}}
\newcommand{\comment}[1]{{\bf COMMENT:\{\{}{#1}{\bf\}\}}}
\newcommand{\plfootnote}[1]{{\footnote{COMMENTAIRE DE PASCAL: {#1}}}}
%
%
\newcommand{\Proof}[1]{{\bf Proof {#1}. }\par}
%
\newcommand{\cqfd} {\mbox{}\nolinebreak\hfill\rule{2mm}{2mm}\medbreak\par}
%
\newcommand{\biindice}[3]%
{
\renewcommand{\arraystretch}{0.5}
\begin{array}[t]{c}
#1\\ {\scriptstyle #2}\\ {\scriptstyle #3}
\end{array}
\renewcommand{\arraystretch}{1}
}
%
\newcommand{\BR}{\Bbb R}         
\newcommand{\BQ}{\Bbb Q}         
\newcommand{\BC}{\Bbb C}         
\newcommand{\BN}{\Bbb N}         
\newcommand{\BZ}{\Bbb Z}         
\newcommand{\BL}{\Bbb L}         
\newcommand{\Bk}{\Bbb k}         
\newcommand{\SSets}{\mathrm{SSets}} 
\newcommand{\Top}{\mathrm{Top}} 
\newcommand{\quism}{\stackrel{\simeq}{\rightarrow}}  
\newcommand{\backquism}{\stackrel{\simeq}{\leftarrow}}
\newcommand{\iso}{\stackrel{\cong}{\rightarrow}}    
\newcommand{\cofarrow}{\to}
\newcommand{\Tor}{\mathrm{Tor}}
\newcommand{\Ext}{\mathrm{Ext}}
\newcommand{\im}{\mathrm{im}}
\newcommand{\pr}{\mathrm{pr}}
\newcommand{\coker}{\mathrm{coker}}
\newcommand{\gldim}{\mathrm{gldim}}
\newcommand{\fdim}{\mathrm{fdim}}
\newcommand{\cat}{{\mathrm{cat}}}
\newcommand{\nil}{{\mathrm{nil}}}
\newcommand{\Cl}{{\mathrm{Cl}}}
\newcommand{\id}{{\mathrm{id}}}
\newcommand{\Sing}{{\mathrm{Sing}}}
\newcommand{\rank}{\mathrm{rank}}
\newcommand{\del}{\partial}

%
\begin{center}%
{\Large{\bf Invariants of t-structures and classification of nullity classes.}
\\ {\normalsize\version}}
\medbreak

{\large\bf Don Stanley}%
\\University of Regina
\end{center}
\vskip5mm
\begin{center}{\bf Abstract}\end{center}
We construct an invariant of t-structures on the derived category of a Noetherian ring.
This invariant is complete when restricting to the category of quasi-coherent complexes,
and also gives a classification of nullity classes with the same restriction.
On the full derived category of $\mathbb Z$ we show that the class of
distinct t-structures do not form a set.
{AMS-classification (2000)} : 18E30, 18E35
\vskip1cm
%
\section{Introduction}
A $t$-structure on a triangulated category generalizes the idea
of truncating the homology of a chain complex above a specified degree.
They were introduced by Beilinson, Bernstein and Deligne in \cite{BBD}.
This paper constructs an invariant of t-structures of
$D(R)$ the derived category of a Noetherian ring $R$.
When restricting attention to the quasi-coherent complexes $D_{\qc}(R)$, we show that
this invariant is complete. We do this by classifying a slightly weaker
structures called a nullity classes.

The Postnikov section functor $P_n$, which kills all homotopy groups
$\pi_i(X)=[S^i,X]$ above dimension $n$, provides a notion of
truncation in topological spaces. Originating in the work of Bousfield \cite{Bousfield1} and
Dror-Farjoun \cite{DF},
for any space $E$ there is a more general truncation functor $P_E$ which kills all maps
from $E$ and it's suspensions so that $[\Sigma^iE, P_E(X)]=0$. It can be defined to
be the universal functor with this property. The class of objects such that
$P_E(X)=0$, $\overline{C(E)}$ is a nullity class. More generally a nullity class is a
full subcategory closed under arbitrary coproducts, (positive) suspensions and extensions.

When working in $D(R)$, in the language of Keller and Vossieck \cite{KV},
a nullity class is a cocomplete preaisle. They showed, more generally
in any triangulated category $\bold T$,
that $t$-structures correspond the preaisles ${\mathcal A}\subset {\bold T}$ allowing a
right adjoint to the inclusion functor. They call such preaisles aisles. Motivated by this and
ideas of Bousfield,  Alonso-Tarrio, Jeremias-Lopez and Souto-Solario \cite{AJS}
constructed the functor $P_E$ in the category $D(R)$.
They also showed that, in our notation, $\overline{C(E)}$ is an aisle. This
gives a very general way to contruct $t$-structures. We will see that when restricted
to quasi-coherent complexes this gives all the nullity classes and thus all the $t$-structures.

Our approach has its roots in topology and thick subcategories. If our nullity class
is also closed under desuspension then it is a localizing subcategory,
an infinite version of a thick subcategory. The thick subcategories of $p$-local finite spectra
were classified by Hopkins and Smith \cite{HS} in terms of an invariant called type. Bousfield
\cite{Bousfield2} used their classification to classify nillity classes of $p$-torsion finite
suspension spaces. Bousfield's classification is in terms of two things: type, which tells us
which thick subcategory the class generates stably
and connectivity, which tells us where the class starts.
Using ideas from the classification for spectra, Hopkins \cite{Hopkins}
and Neeman \cite{Neeman} classified retraction closed
thick subcategories of the perfect complexes
$D_{\parf}(R)$ for Noetherian rings $R$ by subsets of $\Spec R$ closed
under specialization. The invariant is given by taking the support of the object.
Neeman \cite{Neeman} proved the analogous result in $D(R)$ where localizing
subcategories are classified by all subsets of $\Spec R$.

Starting with a nullity class $\mathcal A$ we associate a function (see \ref{p def}):
$$
\phi({\mathcal A})\colon \mathbb{Z}\rightarrow \{ {\rm Subsets\ of\ } \Spec(R)
{\rm \ closed \ under \ specialization }\},
$$
whose value at $n$, $\phi({\mathcal A})(n)$ can be thought of in the following way.
Truncate $\mathcal A$ above $n$ with the standard
truncation to get $\tau^{\geq -n} {\mathcal A}$. Next take the thick subcategory
generated by $\tau^{\geq -n} {\mathcal A}$ and apply the correspondence of Hopkins-Neeman,
in other words take supports. This subset of $\Spec(R)$
is the value $\phi({\mathcal A})(n)$.
Since we cannot desuspend, as in Bousfield's result, we have to prescribe at what
level the $t$-structure starts and there is some choice of when different primes
can start. If $p\in \phi({\mathcal A})(n)$ it means that that prime has already been included at
level $n$. So we see that $\phi({\mathcal A})$ must be increasing.
In this way applying Bousfield's philosophy
to the Hopkins-Neeman result from Theorem \ref{class} we get:

\begin{thmA}
$\phi$ is an order preserving bijection between nullity classes in
$D_{\qc}(R)$ and increasing functions from $\mathbb{Z}$ to subsets of $\Spec R$ closed
under specialization.
\end{thmA}
Since all aisles are nullity classes (Theorem \ref{aislenull}) this theorem also
implies $\phi$ is a complete invariant of $t$-structures in $D_{\qc}(R)$.
Although in $D(R)$ all the $\overline{C(E)}$ are aisles,
their restriction to $D_{\qc}(R)$, $\overline{C(E)}\cap D_{\qc}(R)$ may not
be an aisle. The problem is that for $M\in D_{\qc}(R)$, $P_E(M)$ may no longer
be in $D_{\qc}(R)$. However we get some control over the image, in fact the primes
must be added in a very particular way for the restriction to have a chance
of being an aisle. We get Theorem \ref{comono}:

\begin{thmB}
Suppose $\mathcal A$ is an aisle. Then if $p'\in \phi({\mathcal A})(n)$ and $p$ maximal under $p'$,
then $p\in \phi({\mathcal A})(n+1)$.
\end{thmB}

We conjecture the converse of the theorem: all nullity classes satisfying the
condition are aisles.  Proving the conjecture would complete
the classification of $t$-structures in $D_{\qc}(R)$. It is worth remarking
at this point that the conjecture would be false if we were to work
with perfect complexes (see Example \ref{ex-qcright}), so $D_{\qc}(R)$
seems to be the right place to work when taking this point of view.

In the last section of the paper we use examples of Shelah \cite{Shelah}
to show that there is not a set but rather a proper class of $t$-structures
in $D(\mathbb{Z})$. This not only shows a strong contrast to what happens with
localizing subcategories in $D(R)$ but also shows that classifying the
$t$-structures in $D(R)$ is probably not feasible. These examples
can be transported to the topological setting to show there
exists no set of nullity classes in spectra or topological spaces and
no set of $t$-structures in the triangulated category of spectra.

Next we give a short description of the contents of each section, more details and comments
can be found at the start of some sections. Section \ref{back} gives some background from
ring theory and derived categories and also introduces aisles, nullity classes and
the nullification functor $P_E$. Section \ref{props} proves some properties of nullity classes
and $P_E$ that are well known in the topological setting.
Section \ref{inv} defines the invariant $\phi$ and another $N$ which is the inverse of $\phi$.
In Section \ref{short} we give a short proof that $\phi$ is a complete invariant when restricted
to aisles in $D_{\qc}(R)$. The heart of the paper is Section \ref{sec-null} where the classification
of the nullity classes in $D_{\qc}(R)$ is completed. In Section \ref{image}
we get restrictions on the image
of $\phi$ when restricted to aisles.
Section \ref{propclass} constructs the examples in $D(\mathbb{Z})$ that show the aisles,
and nullity classes do not form a set.

\section{Background and notation}\label{back}
Throughout this paper we let $R$ be a Noetherian ring.
\subsection{Ring theory and associated primes}
\
\vskip 0.3cm
\par
Recall that $\Spec R$ is the set of primes of $R$. If
$U\subset \Spec R$, then $\overline{U}$ denotes the closure of $U$ under specialization.
That is,
$$
\overline{U}=\{ p\in \Spec R | \exists q \in U, q\subset p\}
$$
We will repeatedly use the concept of associated prime. Most things
we need to know about them can be found in Eisenbud's book
\cite{Eisenbud}.
\begin{defin}
For an $R$-module $M$, $\Ass M$ denotes the set of associated primes of $M$.
So $p\in \Ass M$, if $p$ is the annihilator of some element of $M$.
$$
\Ass M=\{ p| \exists  x\in M, \ann x=p\},
$$
where for $x\in M$, $\ann x$ denoted the annihilator of $x$.
\end{defin}

\begin{lemma}\label{WA}
Let $A$, $B$ and $C$ be $R$-modules.
If we have an exact sequence
$$
0 \rightarrow A \rightarrow B \rightarrow C \rightarrow 0
$$
then $\Ass A \subset \Ass B \subset \Ass A \cup \Ass C$.
\end{lemma}
\proof
\cite[Lemma 3.6b]{Eisenbud}.
\cqfd

\begin{lemma}\label{WAW}
Let $B$ and $C$ be $R$-modules.
If we have an exact sequence
$$
B \stackrel{f}{\rightarrow} C \rightarrow 0
$$
then $\Ass C\subset \overline{\Ass B}$.
\end{lemma}
\proof
Let $p\in \Ass C$, and $x\in C$ such that $\ann x=p$.
Let $y\in f^{-1}(x)\subset B$, then clearly $\ann y\subset \ann x=p$. Let
$p'\subset p$ be any prime minimal among primes containing $\ann y$. Then
since the submodule generated by $y$, $\langle y \rangle$ is finitely generated,
$p\in \overline{p'}$ and $p'\in \Ass \langle y \rangle$
by \cite[Lemma 3.1 a]{Eisenbud} and
so $p'\in \Ass B$ by
Lemma \ref{WA}.
Thus $\Ass C\subset \overline{\Ass B}$.
\cqfd

\begin{lemma}\label{A}
Supposing that $M$ is a finitely generated $R$-module,
$p\in \overline{\Ass M}$ if and only if $M\otimes R_p\not=0$.
\end{lemma}
\proof
Assume $p\in \overline{\Ass M}$, then there is a $p'\in \Ass M$ such that $p'\subset p$,
so there is an injective map $R/p'\rightarrow M$. So since $\_\otimes R_{p'}$ is exact
$M\otimes R_{p'}\not=0$ which implies that $M\otimes R_{p}\not=0$.

Next assume $M\otimes R_{p}\not=0$ then there exist $p'\in \Ass (M\otimes R_{p})$ such that
$p'\subset p$ and so $p\in \overline{\Ass M}$ since
$\Ass (M\otimes R_{p})=\{ p'\in \Ass M| p'\subset p\}$ by \cite[Theorem 3.1c]{Eisenbud}
\begin{lemma}\label{B}
Suppose $(R,m)$ is a local ring. If $M$ is a non-trivial finitely generated
$R$-module then there is a surjective map
$M\rightarrow R/m$.
\end{lemma}
\proof
We prove this by quotienting out all but one element of a minimal generating set and then
quotienting out by $m$ times the last generator.
\cqfd
\subsection{Derived categories}
\
\vskip 0.3cm
\par
The category of chain complexes of $R$-modules with differential of
degree $-1$ is denoted by $C(R)$
and
$D(R)$ denotes the derived category of the ring $R$. This is just $C(R)$
modulo weak equivalences. We make $D(R)$ into
a triangulated category in the standard way.
For $M\in D(R)$ or $C(R)$, $H_i(M)$ denotes the $i$-th homology of $M$.
Given $A\in C(R)$, $s^iA$ is the same
complex shifted up by $i$, $s^iA_j=A_{j-i}$ and
$ds^ia=(-1)^is^ida$. Since $s^i(\_)$ preserves weak
equivalences $s^i(\_)$ extends to $D(R)$. By convention we write $s^1=s$.
$H_i(sM)=H_{i-1}(M)$.
If $f\colon A \rightarrow B$ is a map in $D(R)$,
we get a distinguished triangle,
$$
A \rightarrow B \rightarrow C \rightarrow sA.
$$
Applying $H_*$ to a triangle we get a long exact sequence
$$
H_iA \rightarrow H_iB\rightarrow H_iC\rightarrow H_{i-1}A.
$$
If $M$ is an $R$-module we consider it as the object $M$ in $C(R)$ or $D(R)$ with
$$
M_i=
\begin{cases}
M & i=0 \\
0 & else
\end{cases}
$$
and trivial differential.
\par
For a category $\mathcal C$ and $A,B\in {\mathcal C}$
$\Hom_{\mathcal C}(A,B)$ denotes the set of maps from $A$ to $B$. We will omit the subscript
$\mathcal C$ if it is clear which category we are working in, usually it will be $D(R)$.
\par
There are two important full triangulated subcategories of $D(R)$ that we will use:
\par
\noindent
$\bullet \quad D_{\parf}(R)$ is the subcategory of $D(R)$ consisting of objects
represented by chain complexes of finitely generated projective modules.
\par
\noindent
$\bullet \quad D_{\qc}(R)$ is the subcategory of $D(R)$ whose homology
groups are finitely generated and bounded above and below.

\begin{lemma}\label{rag}
If $M,N\in D(R)$, $H_i(M)$ is finitely generated and bounded below, then the natural map
$$
\Hom_{D(R)}(M,N)\otimes R_p\rightarrow \Hom_{D(R_p)}(M\otimes R_p,N\otimes R_p)
$$
is an isomorphism of $R_p$ modules.
\end{lemma}
\proof
Represent $N$ by
$B\in C(R)$ that is bounded above. Since $H_i(M)$ is finitely generated and bounded
below we can represent $M$ by $A\in C(R)$ that is bounded below, projective
in each degree and has finitely many generators in each degree.

Let
$\eHom^i_{grR-mod}(A,B)=\prod_n\Hom_{R-mod}(A_{n+i}, B_n)$,
be the $R$-module of graded $R$-module maps from $A$ to $B$ lowering degree by $i$.
We put a differential on $\eHom_{grR-mod}(A,B)=\oplus \eHom^i_{grR-mod}(A,B)$ by setting
$df(x)=d(f(x))+(-1)^if(dx)$

Then the natural map
$$
\theta\colon \eHom^i_{grR-mod}(A,B)\otimes R_p \rightarrow
\eHom^i_{grR_p-mod}(A\otimes R_p,B\otimes R_p)
$$
is an isomorphism for each $i$ since by \cite[Proposition 2.10]{Eisenbud} it is an
isomorphism restricted to each factor of the product.
The map $\theta$ also commutes with the differential.
Seeing chain maps as $0$-cycles and chain homotopies as $0$-boundaries, we get that
$\Hom_{D(R)}(M,N)=H^0(\eHom_{grR-mod}(A,B))$, and similarly for
$\Hom_{D(R_p)}(M\otimes R_p,N\otimes R_p)$.
Since $R_p$ is flat by \cite[Proposition 2.5]{Eisenbud} for any $L\in C(R)$,
$H_i(L)\otimes R_p\rightarrow H_i(L\otimes R_p)$ is an isomorphism of $R_p$ modules.
So
$$
\Hom_{D(R)}(M,N)\otimes R_p=H^0(\eHom_{grR-mod}(A,B))\otimes R_p
$$
$$
\stackrel{\cong}{\rightarrow} H^0(\eHom_{grR-mod}(A,B)\otimes R_p)
$$
$$
\stackrel{\cong}{\rightarrow} H^0(\eHom_{grR_p-mod}(A\otimes R_p,B\otimes R_p))
$$
$$
=\Hom_{D(R_p)}(M\otimes R_p,N\otimes R_p),
$$
and so the lemma follows.
\cqfd
\begin{lemma}\label{D}
Suppose $M,N \in D_{\qc}(R)$. Suppose
$p\in \overline{\Ass H_n(M)}$ and $p\not\in \overline{\Ass H_i(M)}$ for $i<n$.
If also
$p\in \Ass H_n(N)$, and $p\not\in \Ass H_i(N)$ for $i>n$, then $\Hom(M,N)\not=0$.
\end{lemma}
\proof
Since $M,N \in D_{\qc}(R)$, by Lemma \ref{rag} it is enough to show
that $\Hom_{D(R_p)}(M\otimes R_p, N\otimes R_p)\not=0$. Since
$p\not\in \Ass H_i(M)$ for $i<n$ and $p\in \Ass H_n(M)$,
Lemma \ref{A} implies that $H_i(M\otimes R_p)=0$
for $i<n$ and $H_n(M\otimes R_p)\not=0$. So Lemma \ref{B} implies that there
is a map $f\colon M\otimes R_p \rightarrow s^nR/p$
that induces a surjection on $H_n$. Also
since $p\in \Ass H_n(N)$, \cite[Theorem 3.1 c]{Eisenbud} implies that
$p\in \Ass H_n(N\otimes R_p)$ and thus there is an injection
$g':R/p \rightarrow H_n(N\otimes R_p)$. Since by Lemma
\ref{A} $H_i(N)\otimes R_p=0$ for $i>n$,
there is a map $g: s^nR/p \rightarrow N\otimes R_p$ which induces $g'$
in $H_n$. The composition $gf\colon M\otimes R_p \rightarrow N\otimes R_p$
is nontrivial on $H_n$ and hence nontrivial, thus
$\Hom_{D(R_p)}(M\otimes R_p, N\otimes R_p)\not=0$ and we are done.
\cqfd
Observe that the lemma does need the finiteness condition, as we since since
$\Hom(\mathbb{Q}, \mathbb{Z}/p)=0$ in $D(\mathbb{Z})$. In particular in this case the part of
the proof that relies on Lemma \ref{B} does not hold.
A slight variant of the last lemma is:
\begin{lemma}\label{map}
If $M\in D_{\qc}(R)$ such that $H^n(M)\otimes R_{p}\not=0$ and
$H^i(M)\otimes R_{p}=0$ if $i<n$ then there is a map
$\phi\colon M \rightarrow s^nR/p$ such that $H^n(\phi)\not=0$.
\end{lemma}
\proof
There is a map $M\otimes R_p\rightarrow s^nH^n(M)\otimes R_p$ that is an isomorphism
on $H^n$. Using this map the lemma follows directly from Lemmas \ref{rag} and \ref{B}.
\cqfd
\subsection{Aisles and classes generated by a module}
\
\vskip 0.3cm
\par
Let $\bf T$ be a triangulated category.
We make the give the following definitions due to Keller and Vossieck \cite{KV}.
\begin{defin}\label{tdef}
A full subcategory $\mathcal U$ of $\bf T$
is a {\em pre-aisle} if:
\begin{enumerate}\renewcommand{\labelenumi}{\arabic{enumi})}
\item for every $X\in {\mathcal U}$, $sX\in {\mathcal U}$

\item for every distinguished triangle $X \rightarrow Y\rightarrow Z \rightarrow sX$,
if $X,Z\in {\mathcal U}$ then $Y\in {\mathcal U}$.
\end{enumerate}
A pre-aisle $\mathcal U$ is called {\em cocomplete} if $\mathcal U$ is closed under
coproducts. A pre-aisle $\mathcal U$ such that the inclusion ${\mathcal U}\subset {\bf T}$
admits a right adjoint is called an {\em aisle}.
\end{defin}
Keller and Vosiek \cite{KV} proved that a $t$-structure corresponds to an aisle.
For a definition of $t$-stucture, see \cite{AJS} or \cite{BBD}.
For a subcategory ${\mathcal U}\subset {\bf T}$ we let
$$
{\mathcal U}^{\perp}=\{ x\in {\bf T}|\Hom(y,x)=0 \ \forall y\in {\mathcal U} \}
$$
\begin{thm}\cite{KV}
A pre-aisle $\mathcal U$ is an aisle, that is the inclusion ${\mathcal U}\subset T$
admits a right adjoint, if and only if $({\mathcal U},s{\mathcal U}^{\perp})$ is a
$t$-structure.
\end{thm}

We will mainly consider aisles for the rest of the paper since they
are equivalent to $t$-structures

Cocomplete pre-aisles in the triangulated category of spectra,
and similar subcategories of spaces,
have also been referred to as nullity classes and Bousfield classes.
We will also use the term nullity class to refer to the intersection of a cocomplete
pre-aisle with a full subcategory.

\begin{defin}\label{BC}
Let ${\mathcal D}\subset D(R)$ be a full
triangulated subcategory. A {\em nullity class}
in $\mathcal D$ is a full subcategory of the form ${\mathcal A} \cap {\mathcal D}$
where $\mathcal A$ is
a cocomplete preaisle in $D(R)$. We let $NC$ denote the set of
nullity classes when $\mathcal D=D_{\qc}(R)$.
We order $NC$ by inclusion.
\end{defin}

Notice that the objects in $D(R)$ with finitely generated homology form
a pre-aisle but not a nullity class, so
not all nullity classes are pre-aisles,
however we do have the following:

\begin{thm}\label{aislenull}
Suppose ${\mathcal D}\subset D(R)$ is a full triangulated subcategory.
Any aisle in ${\mathcal D}$ is a nullity class and any nullity class
is a pre-aisle.
\end{thm}

\proof
Since ${\mathcal D}$ is a triangulated subcategory it is clear that
any nullity class is a pre-aisle.
\par
Suppose ${\mathcal U}\subset {\mathcal D}$ is an aisle. By \cite[Proposition 1.1]{AJS}
$x\in {\mathcal U}$ if and only if for every $y\in {\mathcal U}^{\perp}$,
$\Hom(x,y)=0$. This condition is closed under taking coproducts and extensions in the first variable.
Also again by \cite[Proposition 1.1]{AJS}
${\mathcal U}^{\perp}\subset s{\mathcal U}^{\perp}$,
the condition is also closed under suspension. Thus any object in $D(R)$, and
hence in the full subcategory $\mathcal D$, that can be constructed using
these operations also satisfies the condition. The fact that an aisle is a nullity class
follows.
\cqfd
In the proof of the last lemma, the operations can take us out of $\mathcal D$, but that doesn't matter
since we intersect back with $\mathcal D$, and $\mathcal D$ is full.

In \cite{AJS} Alonso Tarrio, Jeremias Lopez and Souto Salorio, show that for any Grothendieck
category $\mathcal A$ and any $E\in D({\mathcal A})$, there is an associated aisle.
A special case of \cite[Proposition 3.2]{AJS} is the following:
\begin{thm}\label{AJS}\cite{AJS}
Let $R$ be a commutative ring and $E\in D(R)$. If $\mathcal U$ is the smallest nullity class
of $D(R)$ that contains $E$, then $\mathcal U$ is an aisle in $D(R)$.
\end{thm}
\vskip 0.5cm
\noindent
{\bf Notation:}
Following the topologists we denote the nullity class
$\mathcal U$ of the proposition associated to $E$
by $\overline{C(E)}$. We also denote the associated truncation functor $\tau_E^{\geq 1}$ by
$P_E$ and $\tau_E^{\leq 0}(X)$ by $X\langle E \rangle$.
More generally for any aisle $\mathcal A$. we will denote
the associated truncation functor $\tau_{\mathcal A}^{\geq 1}$ by
$P_{\mathcal A}$ and $\tau_{\mathcal A}^{\leq 0}(X)$ by $X\langle {\mathcal A} \rangle$.
\par
So for any $X\in D(R)$ we have distinguished triangles
\begin{myeq}\label{long}
X\langle E \rangle \rightarrow X \rightarrow P_E X \rightarrow sX\langle E \rangle.
\end{myeq}
and
$$
X\langle {\mathcal A} \rangle \rightarrow X \rightarrow
P_{\mathcal A} X \rightarrow sX\langle {\mathcal A} \rangle.
$$
The topological notation has the advantage of eliminating the superscripts.
The superscripts are also more compatible with chain complexes with differentials of
degree $+1$, while we are using differentials of degree $-1$. Even though the main reason
I chose this notation is so that I wouldn't get confused, for the purposes of this
paper it seems to be right.

\section{Properties of closed classes and the nullification functor $P_E$}\label{props}.
\begin{lemma}\label{retract}
$\overline{C(E)}$ is closed under retracts.
\end{lemma}
\proof
This follows from the well known
Eilenberg swindle. If $X=A\oplus B$ we can consider the countable coproduct of $X$ with itself
in two different ways.
$\oplus_{i\in \omega} X=A\oplus B\oplus A\oplus B \dots$ or
$\oplus_{i\in \omega} X=B \oplus A\oplus B \oplus A \dots$.
We can include the second into the first missing the first $A$ to get a distinguished
triangle,
$$
\oplus_{i\in \omega} X \rightarrow \oplus X_{i\in \omega} \rightarrow A
\rightarrow s\oplus_{i\in \omega} X
$$
So if $X$ is in $\overline{C(E)}$, since $\overline{C(E)}$ is
cocomplete, so is $\oplus_{i\in \omega} X$ and then
Definition \ref{tdef} 2) implies that $A\in \overline{C(E)}$.
\cqfd
We can put a partial order on $D(R)$ by letting
$E<F \Longleftrightarrow P_E(F)=0$.
The following shows that $<$ is indeed a partial order.
\begin{prop}\label{trans}
$$
E<F \Longleftrightarrow \overline{C(F)}\subset \overline{C(E)}
$$
In particular for any full triangulated subcategory ${\mathcal D}\subset D(R)$,
the map ${\mathcal D}\rightarrow NC$, $E\mapsto \overline{C(E)}$ is order reversing.
\end{prop}
\proof
By definition $E<F$ if and only if $P_EF=0$. Next looking at
the triangle of Equation (\ref{long}), we see that $P_EF=0$ if and only if
$F\langle E\rangle=F$, which happens if and only if
$\overline{C(F)}\subset\overline{C(E)}$.
The second statement follows easily.
\cqfd
As is standard in the topological setting,
we call a map $f$ in $D(R)$ a {\em $P_E$ equivalence} if $P_E(f)$ is an isomorphism in $D(R)$ and
an object $A\in D(R)$, {\em $P_E$ local} if $P_E(A)=A$. The following proposition is standard
in the topological settings and also holds
more generally for any $t$-structure on any triangulated category.

\begin{prop}\label{char uni}
Working in $D(R)$,
\begin{enumerate}\renewcommand{\labelenumi}{\arabic{enumi})}
\item $A \rightarrow P_EA$ is a $P_E$ equivalence.
\item
If $f\colon A \rightarrow B$ is a $P_E$ equivalence then there exists a unique
map $\phi\colon B \rightarrow P_EA$ making the following diagram commute,
$$
\xymatrix{
A \ar[r]^f \ar[dr]_{\eta_A} & B \ar[d]^{\phi} \\
& P_EA.
}
$$
\item $P_EA$ is $P_E$ local.
\item Given $f\colon A \rightarrow B$ with $B$ $P_E$ local, there exists a unique map
$\phi\colon P_EA\rightarrow B$ making the following diagram commute,
$$
\xymatrix{
A \ar[r]^{\eta_A} \ar[dr]_f & P_EA \ar[d]^{\phi} \\
& B.
}
$$
\item $A$ is $P_E$ local if and only if $\Hom(s^iE,A)=0$ for all $i\geq 0$
if and only if $\Hom(X,A)=0$ for all $X\in \overline{C(E)}$.
\item Suppose $E<F$ then $P_E$ local objects are $P_F$ local and
$P_F$ equivalences are $P_E$ equivalences.
\end{enumerate}
\end{prop}
\proof
1):
By Theorem \ref{AJS}, $\overline{C(A)}$ is an aisle in $D(R)$, so it
follows from \cite[Proposition 1.1]{AJS}.
\par\noindent
2):
Since $P_Ef$ is an isomorphism, we can take
$\phi=\eta_B(P_Ef)^{-1}\colon B \rightarrow P_EB\rightarrow P_EA$.
Uniqueness follows from the functoriality of $P_E$.
\par\noindent
3):
Part 1) implies that $P_EA\rightarrow P_EP_EA$ is an equivalence and thus
$P_EA$ is $P_E$ local.
\par\noindent
4):
There is a distinguished triangle
$$
H \rightarrow A \rightarrow P_EA \rightarrow sH
$$
where $H\in \overline{C(E)}$. Since $B$ is $P_E$ local, $\Hom(H,B)=0$ and thus there exists
a dashed extension in the following diagram
$$
\xymatrix{
A \ar[r] \ar[dr]_f & P_EA \ar@{-->}[d]^{\phi} \\
& B.
}
$$
The map $\phi$ is unique since any other extension differs from $\phi$ by an element
of $\Hom(sH,B)$ which is $0$ since $sH\in \overline{C(E)}$.
\par\noindent
5) is \cite[Lemma 3.1]{AJS}.
\par\noindent
6):
For any $A$, we know that $\Hom(s^iE, P_EA)=0$ for every $i\geq 0$, so by Part
5) since $F\in \overline{C(E)}$, $\Hom(s^iF,P_EA)=0$ for every $i\geq 0$, and $P_EA$
is $P_F$ local.
So we get a diagram
$$
\xymatrix
{
A \ar[r]^a \ar[d]_b & P_EA \ar[d]^c \\
P_F A \ar[r]_{P_F(a)} & P_FP_EA
}
$$
in which $c$ is an equivalence. If $A$ is $P_E$ local then $a$ and thus $P_F(a)$ are equivalences.
So by two out of three, $b$ is an equivalence which proves the first part.

For the second part since $P_EA$ is $P_F$ local,
by Part 4) there exists a
dashed extension in the following solid arrow
diagram
$$
\xymatrix
{
A\ar[r]\ar[d] & P_EA \\
P_FA \ar@{-->}[ur]
}
$$
Starting with this diagram and taking $P_E$, in one case more than once,
we get a diagram
$$
\xymatrix
{
A\ar[d] \ar[r] \ar[d] & P_FA \ar[r] \ar[d] & P_EA \ar[d] \\
P_EA \ar[r]_a & P_EP_FA \ar[r]_b & P_EP_EA \ar[r]_c &
P_EP_EP_FA
}
$$
By Part 1) $b\circ a$ and $c\circ b$ are equivalences,
and this implies that $a$ is an equivalence.
\par
Now let $f\colon A\rightarrow B$ be a $P_F$ equivalence.
We get a square
$$
\xymatrix
{
P_EA \ar[r]\ar[d]_{P_Ef} & P_EP_FA \ar[d]^{P_EP_Ff} \\
P_EB \ar[r] & P_EP_FB
}
$$
We have just seen that the horizontal maps are equivalences;
since $f$ is a $P_F$ equivalence $P_Ff$ is an equivalence so $P_EP_Ff$
is an equivalence. So by two out of three $P_Ef$ is an equivalence which is what
we needed to prove.
\cqfd
In $D(R)$ generally direct limits do not exist. Countable homotopy  direct limits in any
triangulated category were constructed in \cite{BN}, and any homotopy direct
limits of objects in  $C(R)$ were constructed in \cite{AJS}. Even though direct limits
in $C(R)$ are homotopy invariant
(this follows since direct limits commute with homology), the direct limits
cannot generally be extended to direct limits in $D(R)$;
phantom maps provide a first obstruction. For these reasons when we do
constructions involving direct limits we will work in $C(R)$.
\begin{prop}\label{closure}
\begin{enumerate}\renewcommand{\labelenumi}{\arabic{enumi})}
\item A direct limit of $P_E$ equivalences is a $P_E$ equivalence.
In particular given a direct system
$\{ A_{\alpha}\}_{\alpha< \lambda}$ of objects in $C(R)$, if
$P_E A(1)\rightarrow P_E A(\alpha)$
is a weak equivalence for each $\alpha<\lambda$ then
$P_E A(1)\rightarrow colim_{\alpha< \lambda}A(\alpha)$ is a weak equivalence.
\item If $A\in \overline{C(E)}$ and
$A\rightarrow B \rightarrow C \rightarrow sA$ is a distinguished triangle then
$B\rightarrow C$ is a $P_E$ equivalence.
\end{enumerate}
\end{prop}
\proof
2): Since $A\in \overline{C(E)}$, $\Hom(A,P_EB)=0$. Thus there exists a dashed
extension $h$ in the following solid arrow diagram
$$
\xymatrix
{
A \ar[d] & \\
B \ar[r] \ar[d]_i & P_EB \ar[d] \\
C \ar[r] \ar@{-->}[ur]^h & P_EC.
}
$$
The map $h$ makes the upper left triangle of the square commute,
and the lower right square commutes too since the two ways around differ
by an element of $\Hom(sA,P_EC)$, which is $0$ since $A\in \overline{C(E)}$.
So we get a commuting diagram
$$
\xymatrix
{
P_E B \ar[r] \ar[d]_{P_Ei} & P_EP_EB \ar[d] \\
P_E C \ar[r] \ar[ur]^{P_Eh} & P_EP_EC
}
$$
in which the horizontal arrow are equivalences.
This implies that $P_Ei$ is an equivalence by Proposition \ref{char uni} 1),
and so $i$ is a $P_E$ equivalence
as required.
\par\noindent
1): In the proof of this part we will be working in $C(R)$.
Consider the direct limit:
$$A(\lambda)=colim_{\alpha<\lambda}A(\alpha).$$
Since $P_E$ is a functor on $C(R)$ we get a commuting diagram,
$$
\xymatrix
{
A(1) \ar[r]\ar[d] & P_E(A_1)\ar[d]\\
A(\lambda)=colim_{\alpha<\lambda}A(\alpha) \ar[r] &
colim_{\alpha<\lambda}P_EA(\alpha).
}
$$
The right vertical arrow is a homology equivalence since it is a
direct limit of homology equivalences.
Thus we get a commuting diagram
$$
\xymatrix
{
A_1 \ar[r] \ar[d] & P_EA_1 \ar[r] & colim_{\alpha<\lambda}P_EA(\alpha) \ar[d] \\
A_{\lambda} \ar[urr] \ar[rr] && P_EA_{\lambda}
}
$$
in which all the horizontal arrows are $P_E$ equivalences. Using Proposition \ref{char uni} 1),
$colim_{\alpha<\lambda}P_EA(\alpha) \rightarrow P_EA_{\lambda}$ is an equivalence
by the same argument as in the proof of Part 2) above. Hence, being a composition of two equivalences,
$P_EA(1) \rightarrow P_E A(\lambda)$ is an equivalence as desired.
\par

\cqfd
The functor $P_E$ can be characterized by its universal properties.

\begin{cor}\label{use}
If $f\colon A\rightarrow B$ is a $P_E$ equivalence and $B$ is $P_E$ local then there is an
isomorphism $\phi\colon B\rightarrow P_EA$ such that
$$
\xymatrix{
A \ar[r]^f \ar[dr]_{\eta_A} & B \ar[d]^{\phi} \\
& P_EA
}
$$
commutes.
\end{cor}
\proof
The map $\phi$ comes from Proposition \ref{char uni} 2), and its inverse from
\ref{char uni} 4). The are compositions are equal to the identity come from Proposition
\ref{char uni} 2) and 4) by the uniqueness part of using universal properties as usual.
\cqfd
\section{An invariant}\label{inv}
We will let $\mathcal S$ denote the set of increasing functions from $\mathbb Z$ to
subsets of $\Spec R$ closed under specialization.
\par

Next we define maps
$$
N\colon {\mathcal S} \rightarrow NC
$$
and
$$
\phi\colon NC \rightarrow {\mathcal S}.
$$

\subsection{Definition of $N$}
Let
$$
{\mathcal S}=\{  f \colon \mathbb Z \rightarrow {\mathcal P}(\Spec R)|
f(n)=\overline{f(n)} \ {\rm and} \ f(n)\subset f(n+1)\},
$$
where $\mathcal P$ is the power set.
We put an order on ${\mathcal S}$ by inclusion, more precisely $f\leq g$ when for every $n$,
$f(n)\subset g(n)$.

For $f\in {\mathcal S}$ let $M(f)=\oplus_n \oplus_{p\in f(n)}s^nR/p$
and $N(f)=\overline{C(M(f))} \cap D_{\qc}(R)$ be the associated nullity class.

\subsection{Definition of $\phi$}\label{p def}
Let ${\mathcal A}\subset D_{\qc}(R)$ be a nullity class.
Define $\phi({\mathcal A})\in {\mathcal S}$ by letting
$p\in \phi({\mathcal A})(n)$ if there is $M\in {\mathcal A}$ such that $p\in \overline{\Ass H^n(M)}$.
So
$$
\phi({\mathcal A})(n)=\{ p\in \Spec R| \exists M\in {\mathcal A} \
{\rm with} \ p\in \overline{\Ass H^n(M)}\}
$$

Note that $\phi({\mathcal A})\in {\mathcal S}$ since if $M\in {\mathcal A}$ and
$p\in \overline{\Ass H^n(M)}$ then $sM\in {\mathcal A}$ and $p\in \overline{\Ass H^{n+1}(sM)}$ also
each $\phi({\mathcal A})(n)$ is clearly closed under specialization from the way they are defined.
Under the correspondence of Hopkins-Neeman \cite{Neeman}
the $\phi({\mathcal A})(n)$ correspond to the thick
subcategories of $D(R)$ generated by the usual truncations of
$\mathcal A$ by dimensions.
\par
As advertised in the abstract $\phi$ can be considered an invariant of $t$-structures in
$D(R)$ by simply intersecting the associated aisle with $D_{\qc}(R)$.
\begin{lemma}\label{compatible}
$N$ and $\phi$ are order preserving.
\end{lemma}
\proof
That $\phi$ is order preserving
follows immediately from the definition. For $f,g \in {\mathcal S}$,
if $f\leq g$ then $M(f)$ is a retract of $M(g)$ so $M(f)\in \overline{C(M(g))}$
by Lemma \ref{retract}. Therefore $\overline{C(M(f))}\subset\overline{C(M(g))}$
and $N$ is seen to be order preserving.
\cqfd

\section{Complete invariant}\label{short}
In this section we give a short proof that when restricted to aisles in $D_{\qc}(R)$,
$\phi$ is injective. This implies that $\phi$ is a complete invariant of such $t$-structures.
We begin with a technical lemma.
\begin{lemma}\label{E}
Let ${\mathcal D}\subset D(R)$ be any full triangulated subcategory.
Let $\mathcal A$ be an aisle in $\mathcal D$ and $M\in {\mathcal D}$. If
$p\in \Ass H_nP_{\mathcal A}(M)$ then $p\in \overline{\Ass H_n(M)}$
or there exists $N\in {\mathcal A}$, $p\in \Ass H_{n-1}(N)$.
\end{lemma}
\proof
There is a distinguished triangle
$$
M\langle {\mathcal A}\rangle \stackrel{f}{\rightarrow} M
\rightarrow P_{\mathcal A}M \rightarrow sM\langle {\mathcal A}\rangle.
$$
From this we get a short exact sequence
$$
0 \rightarrow H_n(M)/\im H(f) \rightarrow H_n(P_{\mathcal A}M)
\rightarrow \ker H_{n-1}(f) \rightarrow 0
$$
Since $\ker H_{n-1}(f)\subset H_{n-1}M\langle {\mathcal A}\rangle$ and
$M\langle {\mathcal A}\rangle\in {\mathcal A}$ ,
the lemma follows directly from Lemmas \ref{WAW} and
\ref{WA}.
\cqfd
\begin{prop}\label{F}
Let ${\mathcal D}=D_{\qc}(R)$ or $D_{\parf}(R)$.
Suppose $\mathcal A$ is an aisle in $\mathcal D$ and $M\in {\mathcal D}$.
Suppose for every $n$ and for every $p\in H_n(M)$, there
exists $l\leq n$, $N\in {\mathcal A}$ and $p\in \overline{\Ass H_l(N)}$, then $M\in {\mathcal A}$.
\end{prop}
\proof
Assume $P_{\mathcal A}M\not=0$. Let $n$ by the largest such that
$H_n(P_{\mathcal A}M)\not=0$. This $n$ exists since $P_{\mathcal A}M\in {\mathcal D}$.
Then there exists a prime $p$ such that
$p\in \Ass H_n(P_{\mathcal A}M)$ and
$p\not\in \Ass H_i(P_{\mathcal A}M)$, if $i>n$.
Thus by Lemma \ref{E}, the hypotheses and the fact
that aisles are closed under suspension, there
exists $N\in {\mathcal A}$ such that $p\in \overline{\Ass H_n(N)}$ and
$p\not\in \overline{\Ass H_l(N)}$ for $l<n$.
So Lemma \ref{D} says that
$\Hom(N,P_{\mathcal A} M)\not=0$,
which contradicts the fact (see \cite[Proposition 1.1]{AJS})
that $P_{\mathcal A}M\in {\mathcal A}^{\perp}$.
So $P_{\mathcal A}M=0$ and $M\in {\mathcal A}$.
\cqfd
It may look like the proof of the last proposition should extend to all nullity classes
or even to any full subcategory ${\mathcal D}\subset D(R)$. Observe though that there are finiteness
conditions needed in the results the proof calls on. In Section \ref{propclass}
we will show there is a proper class of $t$-structures in $D\mathbb{Z}$ so, considering the next
theorem, some finiteness or other assumptions are needed.
\begin{thm}\label{main}
Let ${\mathcal D}=D_{\qc}(R)$ or $D_{\parf}(R)$.
Suppose $\mathcal A, A'$ are aisles in $\mathcal D$, then
${\mathcal A}\subset {\mathcal A}'$ if and only if
for every $n\in \mathbb Z$, $\phi({\mathcal A})(n)\subset \phi({\mathcal A'})(n)$.
Thus $\phi\colon \{{\rm aisles \ in\ } {\mathcal D}\} \rightarrow {\mathcal S}$ is injective.
\end{thm}
\proof
If ${\mathcal A}\subset {\mathcal A}'$ then it is clear from the definition that
$\phi({\mathcal A})(n)\subset\phi({\mathcal A}')(n)$ for all $n$. Also if
$\phi({\mathcal A})(n)\subset\phi({\mathcal A}')(n)$ for all $n$, then it follows directly
from Proposition \ref{F} that ${\mathcal A}\subset {\mathcal A}'$.
\cqfd

\section{Nullity classes in $D_{\qc}(R)$}\label{sec-null}
In this section we classify nullity classes in $D_{\qc}(R)$. This also gives us
another proof that $\phi$ is a complete invariant of $t$-structures.

We use $k(p)$ to denote $(R/p)_{(0)}$.
\begin{lemma}\label{Axfin}
For any finitely generated $R$-module $B$,
$\oplus_{p\in \overline{\Ass B}}R/p< B$.
\end{lemma}
\proof
Since $B$ is finitely generated it has a decomposition
$$
0=B_0\subset B_1\subset \dots \subset B_n=B
$$
such that $B_i/B_{i-1}=\oplus_j R/p(i,j)$ for some primes $p(i,j)$.
By Lemmas \ref{WAW} and
\ref{WA} each $p(i,j)\in \overline{\Ass B}$. Thus Definition
\ref{tdef} 2) then implies that $\oplus_{p\in \Ass B} R/p<B$.
\cqfd
\begin{lemma}\label{Ax}
For any $R$-module $B$,
$\oplus_{p\in \overline{\Ass B}}R/p< B$.
\end{lemma}
\proof
Note that since we will use Proposition \ref{closure} we work in $C(R)$.
Let $\{ x_i \}_{i<\lambda}\subset B$ be a generating set. For $\alpha\leq \lambda$, let
$B(\alpha)\subset B$ be the submodule generated by $\{ x_i \}_{i<\alpha}$. Then
$B=B(\lambda)$ and we will prove the lemma by induction.

Assume $\oplus_{p\in \overline{\Ass B}}R/p< B(\gamma)$ for all $\gamma<\alpha$. If $\alpha$ is
a successor ordinal then $B(\alpha)=colim_{\gamma<\alpha}B(\alpha)$ so
$\oplus_{p\in \overline{\Ass B}}R/p< B(\alpha)$ by Proposition \ref{closure}.

If $\alpha=\gamma+1$ then consider the exact sequence
$$
0 \rightarrow B(\gamma) \rightarrow B(\alpha) \rightarrow M=B(\alpha)/B(\gamma) \rightarrow 0.
$$
The image of $x_{\alpha}$ generates $M$ and so from Lemmas \ref{WAW} and \ref{WA},
$\Ass M\subset \overline{\Ass B(\alpha)} \subset \overline{\Ass B}$.
Thus from Lemma \ref{Axfin}, $\oplus_{p\in \overline{\Ass B}}R/p<M$. So by Definition
\ref{tdef} 2) and the induction hypothesis,
$\oplus_{p\in \overline{\Ass B}}R/p< B(\alpha)$. The lemma now follows by induction.
\cqfd
\begin{lemma}\label{STA}
For every $M\in D(R)$ with homology in only finitely many degrees,
$$
\bigoplus_{i\in {\mathbb Z}}\bigoplus_{p\in\overline{\Ass H^i(M)}} s^iR/p< M
$$
\end{lemma}
\proof
By finiteness $M$ has a decomposition
$$
0\rightarrow M_{r}\rightarrow M_{r+1} \rightarrow \dots \rightarrow M_s=M
$$
such that for every $i$,
$M_i\rightarrow M_{i+1} \rightarrow s^{i+1}H^{i+1}(M)\rightarrow sM_i$
is a distinguished triangle.
Thus it follows easily from Lemma \ref{Ax} and Definition \ref{tdef}, 2) that
$\bigoplus_{i\in {\mathbb Z}}\bigoplus_{p\in\ \overline{\Ass H^i(M)}} s^iR/p< M$.
\cqfd
A similar proof of the last lemma works in the category of bounded above
complexes and presumably the lemma can be proved in the full derived
category using an idea similar to that in Lemma \ref{Ax}
\begin{lemma}\label{torn}
If $E\in D_{\qc}(R)$ then
$$
(P_E M) \otimes R_p\cong P_{E\otimes R_p}(M \otimes R_p)
$$
Where the $P_{E\otimes R_p}$ is taken in the category of $R_p$ modules.
\end{lemma}
\proof
In the construction of $P_EM$ in \cite[Proposition 3.2]{AJS}, there is a cardinal
$\gamma$ and a sequence of objects, $\{ B_{\alpha}\}_{\alpha< \gamma}$ such that
$B(0)=M$, for every $\alpha< \gamma$, there exists a distinguished triangle
$$
\oplus s^kE \rightarrow B_{\alpha} \rightarrow B_{\alpha+1} \rightarrow s\oplus s^kE.
$$
If $\alpha$ is a limit ordinal then $B_{\alpha}=lim_{i<\alpha}B_i$ and
$P_EM=colim_{\alpha< \gamma}B_{\alpha}$. Since $\_\otimes R_p$ preserves triangles
we get a sequence of triangles in $D(R_p)$
$$
\oplus s^kE\otimes R_p \rightarrow B_{\alpha}\otimes R_p \rightarrow B_{\alpha+1}\otimes R_p
\rightarrow s(\oplus s^kE\otimes R_p)
$$
Since $\_\otimes R_p$ commutes with taking colimits
the natural map $colim (B_{\alpha}\otimes R_p)\rightarrow (colim B_{\alpha})\otimes R_p$
is an isomorphism.
These two facts imply that
$M\otimes R_p\rightarrow P_EM\otimes R_p$ is a $P_{E\otimes R_p}$ equivalence.
Also Lemma \ref{rag} implies that $\Hom_{R_p}(E\otimes R_p, P_EM\otimes R_p)=0$.
Thus the lemma follows from Corollary \ref{use}.
\cqfd
\begin{lemma}\label{1A}
$P_AB=0$ implies that for every $M$,
$P_{A\otimes M}B\otimes M=0$.
\end{lemma}
\proof
Recalling the construction of $P_AB$ (see proof of Lemma \ref{torn}) we have a cardinal
$\gamma$, $\{ B_{\alpha} \}_{\alpha< \gamma}$ such that $B(0)=N$ and distinguished triangles
$$
\oplus s^kA \rightarrow B_{\alpha} \rightarrow B_{\alpha+1} \rightarrow s\oplus s^kA.
$$
and if $\alpha$ is a limit ordinal then $B_{\alpha}=colim_{i<\alpha} B_i$
and $P_AB=B_{\gamma}=0$.
Since $\_\otimes M$ preserves triangles and colimits the result follows.
\cqfd
\begin{lemma}\label{1B}
For $M\in D(R)$, if $H_*(M)=0$ for $*<0$, then
$P_RM=0$.
\end{lemma}
\proof
Straightforward.
\cqfd
\begin{lemma}\label{good1}
If $P_AB=0$ and $H_i(M)=0$ for $i<0$ then
$P_A(B\otimes M)=0$.
\end{lemma}
\proof
Lemma \ref{1B} says that $R<M$ so Lemma \ref{1A} implies that $B<B\otimes M$. Since
$A<B$ by assumption, the transitivity of $<$ (Proposition \ref{trans})
implies that $A<B\otimes M$ and we are done.
\cqfd
\begin{lemma}\label{1C}
If $M\in D_{\qc}(R)$ and $q\in \overline{\Ass H_0(M)}$, then
$M<k(q)$.
\end{lemma}
\proof
By Lemma \ref{1B}, $R<k(q)$. So by Lemma \ref{1A}
$M<M \otimes k(q)$. By \cite[Lemma 2.17]{BN},
$M \otimes k(q)$ is a direct sum of suspensions of
$k(q)$. In degree $0$ this direct sum is non-trivial
by Lemma \ref{A}, since $H_0(M)$ is finitely generated and
$q\in \overline{\Ass H_0(M)}$. The result follows from Lemma \ref{retract}.
\cqfd
The last lemma does not always hold for $M\in D(R)$ as we see by taking
$R=\mathbb{Z}$, $M=\mathbb{Q}$ and $q=(p)$ for any non zero prime $p\in \mathbb{Z}$.
\begin{lemma}\label{LL}
$\Ass (k(q))=\{ q \}$.
\end{lemma}
\proof
Let $\frac{x}{u}\in k(q)$ with $x\in R/q$ and $u\in R\setminus q$.
Clearly $p\subset \ann x$ and if $l\in \ann x$ then $vlx=0$ for some $v\in R\setminus q$.
Since $R/q$ is an integral domain this implies either $x=0$ or $l\in q$. So $\ann x=R$ or
$\ann x\subset p$ and we are done.
\cqfd

\begin{prop}\label{1D}
Suppose $\dim(R)$ is finite and  $S\in D_{\qc}(R)$.
For every $p\in \Ass H_0(S)$ and $q$ such that $p\subset q$, $S<s^{\dim R/q}R/q$.
In particular $S<s^{\dim R}R/q$.
\end{prop}
\proof
Fix $p\in \Ass H_0(S)$. Looking at a particular $q$,
assume for every $q'$ with $q\subset q'$, $q'\not=q$ the lemma holds.
Let $M$ be defined to make the following sequence short exact
$$
0\rightarrow R/q \rightarrow k(q) \rightarrow M \rightarrow 0
$$
By Lemma \ref{LL}, $\Ass (k(q))=\{ q \}$ and so by Lemma \ref{WAW},
$\Ass M \subset \overline{\Ass (k(q))}=\overline{q}$.
Since $R/q\otimes k(q) \rightarrow k(q)\otimes k(q)$ is
an isomorphism, $M\otimes k(q)=0$.
If $q\in \Ass M$ then there exists an injection $R/q\rightarrow M$ and this would
mean that $M\otimes k(q)\not=0$. So $q\not\in \Ass M$, and
$\Ass M\subset \overline{q}-\{ q\}$. Therefore
by the induction hypothesis
and Lemma \ref{Ax} $S<s^{\dim R/q-1}M$.
By Definition \ref{tdef} 2) and Lemma \ref{1C},
$S<k(q)<s^{\dim R/q-1}k(q)$ so by the short exact
sequence above $S<s^{\dim R/q}R/q$.
Notice that if $\dim R/q=0$ then $q$ is maximal and
$k(q)=R/q$, so $S<k(q)=R/q$. This completes the proof of the first
statement of the proposition. The second statement follows since
$\dim R/q \leq \dim R/p$.
\cqfd
\begin{lemma}\label{small ring}
Suppose that $\dim R$ is finite and
$M\in D_{\qc}(R)$. If $p'\in ass H^n(M)$ and $p'\subset p$ then
there exist $N\in D_{\qc}(R)$ such that $p\in \Ass H^n(N)$, $M<N$ and
for every $i$, if $p''\in \Ass H^i(N)$ then $p\subset p''$.
\end{lemma}
\proof
Since $M<sM$ we can assume that $n$ is the smallest such that there is a
$p'\in \Ass H^n(M)$ with $p'\in p$. Thus $H^i(M\otimes R_p)=0$ for
$i<n$. Suppose $p=(x_1, \dots, x_s)$ and let
$K=K(x_1, \dots, x_s)=\otimes_i {\rm Cone}(R\stackrel{\cdot x_i}{\rightarrow} R)$
denote the Koszul complex. Let $N=M\otimes K$. By Lemma \ref{good1} $M<N$.
By \cite[Proposition 17.14]{Eisenbud}
if $y\in p$ then $y$ annihilates $H^*(N)$, hence
for every $i$, if $p''\in \Ass H^i(N)$ then $p\subset p''$. Using this and the
fact that
$H^i(M\otimes R_p)=0$ for $i<n$ we calculate that $p\in \Ass H^n(N)$. This shows that
$N$ satisfies the desired conditions.
\cqfd

\begin{lemma}\label{LM1}
Suppose $\dim(R)$ is finite and $M\in D_{\qc}(R)$. If $p\in \Ass H^i(M)$ and
$p\subset p'$ then $M< s^tR/p'$ for all $t\geq i$.
\end{lemma}
\proof Fix a prime $p$ and assume that the lemma is true for
each $M\in D_{\qc}(R)$ and each prime $p''$ such that $p\varsubsetneq
p''$. We wish to show the lemma is true for $p$. So assume
that $p\in \Ass H^i(M)$. By Lemma \ref{small ring} and the
induction hypothesis if $p\varsubsetneq p'$ then $M< s^iR/p'$. So
we just have to prove that $M< s^iR/p$.
\par
From Proposition \ref{1D} we know that for some $k$,  $M< s^kR/p$.
Let $j\leq k$ be the smallest
such that $M< s^jR/p$. If $j\leq i$ we are done, otherwise all that
remains is to show that $M<s^{j-1}R/p$. Since $M<sM$ and $<$ is
transitive using Lemma \ref{small ring}, there exists $N$ such
that $M<N$, $H^{j-1}(N)\otimes R_p\not=0$,$H^l(N)\otimes R_p=0$ if
$l<j-1$ and if $p''\in \Ass H^n(N)$ for some $n$ then $p\subset
p''$.

By Lemma \ref{map} there is a map $\phi\colon N \rightarrow
s^{j-1}R/p$ such that $H^{j-1}(\phi)\not=0$. A simple calculation
with the long exact sequence on homology then shows that
$\Ass H^t(C(\phi))\subset \overline{\Ass H^t(N)\cup \Ass H^{t-1}(N)}$
and $\Ass H^{j-1}(C(\phi))\subset \overline p -p$.
 Thus the induction hypothesis says that for
 each $t$ and $p\in \Ass H^t(C(\phi))$, $M<s^tR/p$.
Thus by Lemma \ref{STA} we get that $M<C(\phi)$, since also $M<N$
we get that $M<s^{j-1}R/p$ and we are done.
\cqfd
Next we remove the hypothesis that $\dim R$ is finite by reducing to the local case and using the last
lemma.
\begin{lemma}\label{LM}
Suppose $M\in D_{\qc}(R)$. If $p\in \Ass H^i(M)$ and
$p\subset p'$ then $M< s^tR/p'$ for all $t\geq i$.
\end{lemma}
\proof
Let $q$ be any prime of $R$ then by Lemma \ref{torn}
$$
P_M(s^tR/p')\otimes R_q \cong P_{M\otimes R_q}s^t(R/p'\otimes R_q).
$$
If $p'\subset q$ then $p\subset q$ and $p\otimes R_q\in \Ass H^i(M\otimes R_q)$.
Also $p\otimes R_q \subset p'\otimes R_q$ and $(R/p')\otimes R_q=R/(p'\otimes R_q)$.
By the Krull Principal Ideal Theorem (\cite[Theorem 10.2]{Eisenbud}) $\dim R_q$ is finite, so by
Lemma \ref{LM1} $P_{M\otimes R_q}(s^tR/q'\otimes R_q)=0$. So we have that for all primes
$q$ of $R$, $P_M(s^tR/p')\otimes R_q=0$ and therefore by \cite[Lemma 2.8]{Eisenbud},
$P_M(s^tR/p')=0$. By definition this is the same as saying that
$M<s^tR/p'$.
\cqfd
\begin{prop}\label{LA}
For any nullity class ${\mathcal A}$ with homology in finitely many degrees,
${\mathcal A}\subset N\phi({\mathcal A})$.
\end{prop}
\proof
Let $M\in {\mathcal A}$. Then by definition for every $n\in \mathbb{Z}$ and
$p\in \overline{\Ass H^n(M)}$, $p\in \phi({\mathcal A})(n)$. Thus
$M(\phi({\mathcal  A}))<\oplus_n\oplus_{p\in \overline{\Ass H^n(M)}}R/p$ since
it is a retract of $M(\phi({\mathcal A}))$. So Lemma \ref{STA} and Proposition \ref{trans}
imply that
$M(\phi({\mathcal  A}))<M$ and therefore $M\in N\phi{\mathcal A}$, so
${\mathcal A}\subset N\phi({\mathcal A})$.
\cqfd
\begin{prop}\label{LB}
For every nullity class ${\mathcal A}\subset D_{\qc}(R)$,
$N\phi({\mathcal A})\subset {\mathcal A}$.
\end{prop}
\proof
Suppose $p\in \phi({\mathcal A})(n)$, then there exist $M(p,n)\in {\mathcal A}$
with $p\in \overline{\Ass H^n(M)}$.
By Lemma \ref{LM}, $M(p,n)<s^nR/p$. Therefore
$\oplus_n\oplus_{p\in \phi({\mathcal A})(n)} M(p,n) \in {\mathcal A}$ and
$$
\oplus_n\oplus_{p\in \phi({\mathcal A})(n)} M(p,n)<
\oplus_n\oplus_{p\in \phi({\mathcal A})(n)} s^nR/p=M(\phi({\mathcal A})).
$$
So $M(\phi({\mathcal A}))\in
\overline{C(\oplus_n\oplus_{p\in \phi({\mathcal A})(n)} M(p,n))}\subset {\mathcal A}$.
and $N\phi({\mathcal A})\subset {\mathcal A}$ as desired.
\cqfd
Notice that the condition ${\mathcal A}\subset D_{\qc}(R)$ is needed, since
$(0)\in \Ass \mathbb{Q}$, so $\mathbb{Z}\in N\phi(\overline{C(\mathbb{Q})})$.
However $\mathbb{Q}\not<\mathbb{Z}$ and so $\mathbb{Z}\not\in \overline{C(\mathbb{Q})}$.
\begin{prop}\label{LC}
Working in $D_{\qc}(R)$, for any $f\in {\mathcal S}$, $\phi Nf=f$.
\end{prop}
\proof
Suppose $p\in f(n)$ then $s^nR/p\in Nf$ and so since
$p\in \Ass H^n(s^nR/p)$, $p\in \phi Nf(n)$.
\par
Now suppose $p\in \phi Nf(n)$. Then there is a $M\in N(f)$ such that
$M(f)<M$ and $p\in {\overline \Ass H^n(M)}$. Thus Lemma \ref{A} and
\cite[Lemma 2.17]{BN}, which says that
$M\otimes k(p)$ is a direct sum of suspensions of $k(p)$, imply that
$s^nk(p)$ is a retract of $M\otimes k(p)$. So using Lemmas \ref{retract} and \ref{1A}
and Proposition \ref{trans}, we see that
$M(f)\otimes k(p)<M\otimes k(p)<s^nk(p)$.
Since $M(f)\otimes k(p)$ is also a direct sum of suspensions of
$k(p)$, it follows that for some $l\leq n$,
$s^lk(p)$ is a retract of $M(f)\otimes k(p)$. Next we applying Lemma \ref{A},
$p\in \overline{H^lM(f)}$ and so $p\in f(l)$ since $f(l)$ is closed
under specialization. Since $f$ is increasing and $l\leq n$,
$p\in f(n)$.
\cqfd
The next theorem provides a classification of nullity classes in $D_{\qc}(R)$.
\begin{thm}\label{class}
\par
$\phi\colon NC\rightarrow {\mathcal S} $ and
$N\colon {\mathcal S} \rightarrow NC $ are inverse bijections of partially ordered sets.
\end{thm}
\proof
This follows from the last three lemmas \ref{LA}, \ref{LB} and \ref{LC}.
\cqfd
As a corollary we give another proof of Theorem \ref{main}
\begin{cor}\label{cor-inj}
Suppose $\mathcal A, A'$ are aisles in $D_{\qc}(R)$, then
${\mathcal A}\subset {\mathcal A}'$ if and only if
for every $n\in \mathbb Z$, $\phi({\mathcal A})(n)\subset \phi({\mathcal A'})(n)$.
Thus $\phi\colon aisles \rightarrow {\mathcal S}$ is injective.
\end{cor}
\proof
Follows from Theorem \ref{class}.
\cqfd
We can consider constant functions $f\in {\mathcal S}$ so that $f(i)=f(j)$ for
all $i,j\in \mathbb{Z}$. Taking $N$ of such a function we get a nullity class
$N(f)$ that is closed under desuspension. As such it is a thick subcategory in
$D_{\qc}(R)$ that is closed under retracts, yet we do not get all
retraction closed thick
subcategories in this way. For example consider $R=\mathbb{Z}/4$.
$\Spec \mathbb{Z}/4=\{ (2) \}$, so a constant function
$f\colon \mathbb{Z} \rightarrow \Spec \mathbb{Z}/4$ is either $f(n)=\emptyset$ or
$f(n)=\{ (2) \}$. If $f(n)=\emptyset$ then $N(f)=\{ 0 \}$, the
class consisting of only the trivial complex, and if $f(n)=\{ (2) \}$ then
$N(f)=D_{\qc}(R)$. However $D_{\parf}(R)\subset D_{\qc}(R)$ is another thick subcategory
that is closed under retracts, and $D_{\parf}(R)\not= D_{\qc}(R)$ since
$\mathbb{Z}/2\notin D_{\parf}(R)$ as it only has infinite resolutions. So in
$D_{\qc}(\mathbb{Z}/4)$ there are more thick subcategories closed under retracts than
nullity classes close under suspension. Nevertheless considering
constant functions in $\mathcal S$ simply as subsets of $\Spec R$ we do get the following corollary of
Theorem \ref{class}.
\begin{cor}
$\phi$ and $N$ induce order preserving bijections between the set of nullity classes
in $D_{\qc}$ closed under desuspension and subsets of $\Spec R$ closed under specialization.
\end{cor}
\proof
Follows directly from Theorem \ref{class} and the definitions of $\phi$ and $N$.
\cqfd
It is tempting to think that by restricting to $D_{\parf}(R)$, we should be able to recover
the result of Hopkins and Neeman, but I know of no way to do that.

\section{Image of invariant}\label{image}
In this section using the classification of nullity classes
we get some control over what the image of $\phi$ is when restricted
to $t$-structures.
The main object of the section is to show that if ${\mathcal A}\subset D_{\qc}(R)$ is an aisle then
then if $p\in \phi({\mathcal A})(n)$ then all primes maximal under $p$ must
be in $\phi({\mathcal A})(n+1)$ (Theorem \ref{comono}). So $\phi({\mathcal A})$
must increase in a very particular way.
\par
As a motivating example let us work in $D(\mathbb{Z}_{(p)})$ and consider
$P_{\mathbb{Z}/p}s\mathbb{Z}_{(p)}$. Since we have a short exact sequence
$$
0\rightarrow \mathbb{Z}_{(p)} \stackrel{\times p}{\rightarrow} \mathbb{Z}_{(p)}
\rightarrow \mathbb{Z}/p \rightarrow 0
$$
we have a triangle
$$
\mathbb{Z}/p \rightarrow s\mathbb{Z}_{(p)} \stackrel{\times p}{\rightarrow} s\mathbb{Z}_{(p)}
\rightarrow s\mathbb{Z}/p
$$
so $s\mathbb{Z}_{(p)} \stackrel{\times p}{\rightarrow} s\mathbb{Z}_{(p)}$ is a
$P_{\mathbb{Z}/p}$ equivalence. Taking colimits we see that
$$
s\mathbb{Z}_{(p)} \rightarrow colim(s\mathbb{Z}_{(p)}
\stackrel{\times p}{\rightarrow} s\mathbb{Z}_{(p)}
\stackrel{\times p}{\rightarrow} \dots)= s\mathbb{Q}
$$
is a $P_{\mathbb{Z}/p}$ equivalence. Also $\Hom(s^i\mathbb{Z}/p,s\mathbb{Q})=0$ for all $i\geq 0$, so
$P_{\mathbb{Z}/p}s\mathbb{Z}_{(p)}= s\mathbb{Q}$. However
$s\mathbb{Z}_{(p)}\in D_{\qc}(\mathbb{Z}_{(p)})$ but
$s\mathbb{Q}\not\in D_{\qc}(\mathbb{Z}_{(p)})$. This implies that the nullity class
${\mathcal A}=\overline{C(\mathbb{Z}/p)}$ is not an aisle in $D_{\qc}(\mathbb{Z}_{(p)})$.
In fact we would need $s\mathbb{Z}_{(p)}\in {\mathcal A}$ to make it an aisle.
It is this basic phenomenon that stops many nullity classes from being aisles.

Recall from \ref{p def} that for $f \in {\mathcal S}$,
$M(f)=\oplus_i\oplus_{p\in f(i)}s^iR/p$.
\begin{lemma}\label{Ming}
Suppose $R$  is a local ring with maximal ideal $m$, and $p$ a prime maximal under $m$.
Let $h\in m-p$ be any element.
Suppose $f\in {\mathcal S}$ such that $m\in f(n)$, $f(n-1)=\emptyset$ and
$p\not\in f(n+1)$.
If $N=s^n(R/p)/(h)\oplus_{q\not=p\in f(n)}s^nR/q
\bigoplus\oplus_{i\not=n}\oplus_{q\in f(i)}s^iR/q$
then $P_{M(f)}=P_N$.
\end{lemma}
Notice to get $N$ from $M(f)$ we simply replaced $s^nR/m$ by $s^n(R/p)/(h)$ and left
everything else the same.
\par
\noindent
\Proof{of lemma}
The only prime that contains $(h)=\ann((R/p)/(h))$ is $m$ and therefore by
\cite[Theorem 3.1 a)]{Eisenbud}, $\Ass((R/p)/(h))=m$. Therefore by Theorem \ref{class},
$s^n(R/p)/(h)<s^nR/m$ and $s^nR/m<s^n(R/p)/(h)$. It follows that $M(f)<N$ and $N<M(f)$, therefore
$P_{M(f)}=P_N$.
\cqfd
\begin{lemma}\label{Ping}
Suppose $R$  is a local ring with maximal ideal $m$, and $p$ a prime maximal under $m$.
Let $h\in m-p$ be any element.
Suppose $f\in {\mathcal S}$ such that $m\in f(n)$, $f(n-1)=\emptyset$ and
$p\not\in f(n+1)$.
$P_{M(f)}(s^{n+1}R/p)=s^{n+1}(R/p[\frac{1}{h}])$
\end{lemma}
\proof
$P_{M(f)}=P_N$ from Lemma \ref{Ming}.
Suppose $q\in f(n)$ or $f(n+1)$. If $q\subset p$ then since $f(n)$ and $f(n+1)$ are
closed under specialization and $f(n)\subset f(n+1)$, we see that $p\in f(n+1)$.
\par
So $q\not\subset p$ and we can choose $h\in q$. So for any map
$f\colon s^iR/q\rightarrow s^{n+1}R/p[\frac{1}{h}]$ we get a square,
$$
\xymatrix
{
s^i R/q \ar[r]^f \ar[d]_{\times h} & s^{n+1}R/p[\frac{1}{h}] \ar[d]^{\times h} \\
s^i R/q \ar[r]^f  & s^{n+1}R/p[\frac{1}{h}].
}
$$
The left vertical map is $0$ but the right vertical map is an isomorphism, this implies
that $f=0$. So $\Hom(s^iR/q, s^{n+1}R/p[\frac{1}{h}])=0$.
It follows that $s^{n+1}(R/p[\frac{1}{h}])$ is $M(f)$ local.
Considering the exact sequence
$$
0 \rightarrow R/p \stackrel{\times h}{\rightarrow} R/p \rightarrow (R/p)/(h) \rightarrow 0
$$
We get a triangle
$$
s^n(R/p)/(h) \rightarrow s^{n+1}R/p \stackrel{\times h}{\rightarrow} s^{n+1}R/p \rightarrow s(R/p)/(h)
$$
And so Proposition \ref{closure} 2) implies that
$s^{n+1}R/p \stackrel{\times h}{\rightarrow} s^{n+1}R/p$ is a $P_{s^n(R/p)/(h)}$
equivalence. Since $M(f)<N<P_{s^n(R/p)/(h)}$ it is also a
$P_{M(f)}$ equivalence by Proposition \ref{char uni} 6).
Since $s^{n+1}R/p[\frac{1}{h}]$ is the colimit of such maps,
it follows from
Proposition \ref{closure} 1) that
$s^{n+1}R/p\rightarrow s^{n+1}(R/p[\frac{1}{h}])$
is a $P_{M(f)}$ equivalence.
Thus, since we saw above that $s^{n+1}(R/p[\frac{1}{h}])$ is $M(f)$ local
$P_{M(f)}(s^{n+1}R/p)=(s^{n+1}R/p[\frac{1}{h}])$ by Corollary \ref{use}.
\cqfd
\begin{lemma}\label{Ling}
Using notation from last few lemmas,
$s^{n+1}(R/p[\frac{1}{h}])\not\in D_{\qc}(R)$.
\end{lemma}
\proof
We know that
$R/p[\frac{1}{h}]=colim (R/p \stackrel{\times h}{\rightarrow}
R/p \stackrel{\times h}{\rightarrow} \cdots )$.
Since $R/p$ is an integral domain, each map $\times h$ is injective and since
$h\in m\setminus p$, $h$ is not a unit and so $\times h$ is not surjective.
These two facts imply that $R/p[\frac{1}{h}]$ is not a finitely generated $R$-module.
Thus since $H^{n+1}(s^{n+1}R/p[\frac{1}{h}])=R/p[\frac{1}{h}]$,
$s^{n+1}R/p[\frac{1}{h}]\not\in D_{\qc}(R)$.
\cqfd
\begin{prop}\label{not finite}
Suppose $p'\in \Spec R$ and $p$ a prime maximal under $p'$.
Suppose $f\in {\mathcal S}$ such that $p'\in f(n)$ and
$p\not\in f(n+1)$, then $P_{M(f)}(s^{n+1}R/p)\not\in D_{\qc}(R)$.
\end{prop}
\proof
By Lemmas \ref{torn} and \ref{Ping}
$$
(P_{M(f)}s^{n+1}R/p\otimes R_{p'}\cong
P_{M(f)\otimes R_{p'}}(s^{n+1}R/p\otimes R_{p'})\cong
s^{n+1}R/p\otimes R_{p'}[\frac{1}{h}]
$$
By Lemma \ref{Ling}, $s^{n+1}R/p\otimes R_{p'}[\frac{1}{h}]\not\in D_{\qc}(R_{p'})$.
Hence $P_{M(f)}s^{n+1}R/p\otimes R_{p'}\not\in D_{\qc}(R_{p'})$ and so
$P_{M(f)}s^{n+1}R/p\not\in D_{\qc}(R)$.
\cqfd
\begin{thm}\label{goodfinite}
Suppose $E\in D(R)$ and ${\mathcal D}\subset D(R)$ is a full triangulated subcategory.
If for every $M\in {\mathcal D}$, $P_EM\in {\mathcal D}$ then the nullity class
${\mathcal A}=\overline{C(E)}\cap {\mathcal D}$ is an aisle.
The converse is also true if there exist a set
$\{ E(\alpha)\}_{\alpha<\lambda}$ of objects in
$\mathcal A$ such that $\oplus_{\alpha<\lambda}E(\alpha)<E$.
\end{thm}
\proof
Suppose that for every $M\in {\mathcal D}$, $P_EM\in {\mathcal D}$. Looking
at the distinguished triangle of Equation \ref{long},
 we see that $M\langle E \rangle \in {\mathcal D}$ as well.
The functor $M\mapsto M\langle E \rangle$ gives the required right adjoint to
the inclusion ${\mathcal A}\subset {\mathcal D}$ and so ${\mathcal A}$ is an aisle.
\par
Now suppose $\mathcal A$ is an aisle. Let $M\in {\mathcal D}$. We know,
see \cite[Proposition 1.1]{AJS} for example for a proof, that we have a triangle
in $\mathcal D$
$$
M\langle {\mathcal A} \rangle \rightarrow M \rightarrow P_{\mathcal A} M \rightarrow
sM\langle {\mathcal A} \rangle
$$
such that:
\par
a) $M\langle {\mathcal A} \rangle\in {\mathcal A}$.
\par
b) $P_{\mathcal A}M\in {\mathcal A}^{\perp}$.
\par
By Proposition \ref{closure} 2), a) implies that $M \rightarrow P_{\mathcal A} M$
is a $P_{E}$ equivalence.
\par
Statement b) above says that for every $X\in {\mathcal A}$, $\Hom(X,P_{\mathcal A}M)=0$.
In particular for every $\alpha$, since $E(\alpha)\in {\mathcal D}$
and $E<E(\alpha)$, $E(\alpha)\in {\mathcal A}$, and so
$\Hom(E(\alpha),P_{{\mathcal A}}M)=0$. Thus
$$
\Hom(\oplus_{\alpha<\lambda}E(\alpha),P_{\mathcal A}M)=
\prod_{\alpha<\lambda}\Hom(E(\alpha),P_{\mathcal A}M)=0,
$$
and so $P_{\mathcal A}M$ is $P_{\oplus_{\alpha<\lambda}E(\alpha)}$ acyclic
and thus since $\oplus_{\alpha<\lambda}E(\alpha)<E$,
$P_E$ acyclic by Proposition \ref{char uni} 6).
So by Corollary \ref{use},
$P_{E}M\cong P_{\mathcal A}M\in {\mathcal D}(R)$.
\cqfd
The condition that $\oplus_{\alpha<\lambda}E(\alpha)<E$ arises
since something could be in $\overline{C(E)}^{\perp}$
when restricted to maps in a smaller category, like $\mathcal D$, but no longer in
$\overline{C(E)}^{\perp}$ in $D(R)$. This seems related to the problem of the construction
of cohomological localizations in the category of spectra.
\begin{cor}\label{finite}
For $f\in {\mathcal S}$, $N(f)$ is an aisle if and only if for every
$A\in D_{\qc}(R)$, $P_{M(f)}A\in D_{\qc}(R)$.
\end{cor}
\proof
By definition $N(f)=\overline{C(M(f)}\cap D_{\qc}(R)$ and
$M(f)=\oplus_n\oplus_{p\in f(n)} s^nR/p$, so in particular
$\oplus_n\oplus_{p\in f(n)} s^nR/p<M(f)$ and each $s^nR/p\in N(f)$. Thus the corollary follows
from the theorem.
\cqfd
We call a function $f\colon \mathbb{Z}\rightarrow \Spec R$ {\em comonotone} if whenever
$p'\in f(n)$ and $p$ is maximal under $p'$, then $p\in f(n+1)$.
\begin{thm}\label{comono}
If $\mathcal A$ is an aisle in $D_{\qc}(R)$, then $\phi({\mathcal A})$ is comonotone.
\end{thm}
\proof
Follows directly from Proposition \ref{not finite} and Theorem \ref{finite}.
\begin{conj}
For a noetherian ring $R$, if $f$ is comonotone then $N(f)$ is an aisle.
\end{conj}
The converse of this is Theorem \ref{comono} and by Corollary \ref{cor-inj}
or Theorem \ref{main}, all
aisles in $D_{\qc}(R)$ are of this form. So proving this conjecture would complete the classification
of $t$-structures in $D_{\qc}(R)$.
\vskip 0.3cm
\begin{example}\label{ex-qcright}
Quasi-coherent complexes are really the right ones to work
with in this case since there are additional
restrictions for having a $t$-structure in $D_{\parf}(R)$. For example
consider again $D_{\parf}(\mathbb{Z}/4)$ and
$f\colon \mathbb{Z} \rightarrow \Spec \mathbb{Z}/4=\{ (2)\}$ given by
$$
f(i)=
\begin{cases}
\emptyset & i\leq 0 \\
\{ (2) \} & i>0
\end{cases}
$$
Then $M(f)=\oplus_{i>0}s\mathbb{Z}/2$, $N(f)$ is just all complexes with homology concentrated
in possitive degrees and $P_{M(f)}$ is just truncation. Letting $A$ be the complex
$$
A_i=
\begin{cases}
\mathbb{Z}/4 & i=0,1 \\
0 & {\rm else}
\end{cases}
$$
and $d\colon A_1\rightarrow A_0$ be multiplication by $2$,
we can see that $P_{M(f)}A=\mathbb{Z}/2$, but
$\mathbb{Z}/2\not\in D_{\parf}(\mathbb{Z}/4)$ since any resolution of it has infinite length.
Also $sA\in \overline{C(M(f)} \cap D_{\parf}(R)$, and $sA<s\mathbb{Z}/2$, so we get by
Theorem \ref{goodfinite} that $\overline{C(M(f))}\cap D_{\parf}(R)$ is not an aisle.
\end{example}

\section{A class of t-structures in $D(\mathbb Z)$.}\label{propclass}

In this section we show that the t-structures in $D(\mathbb Z)$ do not form a
set but rather a proper class (Corollary \ref{properclass}). The same proof shows that the nullity classes
in spectra and in topological spaces do not form a set. Similarly the $t$-structures in the triangulated
category of spectra do not form a set.
These results follows easily from some nice,
and more difficult, examples of Shelah \cite{Shelah}.
\par
There are two main reasons we chose to exhibit this result. The first reason is to show that it
is unreasonable to expect a nice classification of $t$-structures or nullity classes in $D(R)$.
The second, related, reason is to contrast with what happens in the case of
localizing subcategories in $D(R)$.
If we demand that our nullity classes are also closed under taking desuspensions, we get a
localizing category in $D(R)$. Neeman \cite {Neeman} showed that these are in 1-1 correspondence
with subsets of $\Spec R$, so the situation is only slightly more complicated than for
thick subcategories of $D_{\parf}(R)$. So going from something with
some finiteness conditions, $D_{\parf}(R)$, to infinite things, $D(R)$, only increases
complexity slightly. However the situation for nullity classes is much
different. In the finite case, $D_{\qc}(R)$, we have a classification more or less
in terms of increasing sequences of thick subcategories; when we move to $D(R)$ though, we
completely lose control, we might have a proper class of nullity classes, and even with the extra
condition required for a $t$-structure still have a proper class.
\begin{defin}
A system $\{ A_{\alpha}\}_{\alpha \in Y}$ of distinct abelian groups is called a rigid system
if $\alpha\not=\beta$ implies that
$\Hom(A_{\alpha},A_{\beta})=0$.
\end{defin}
In \cite{Shelah}, Shelah proved the following:
\begin{thm}\label{She}\cite{Shelah}
For every cardinal $\lambda$ there is a rigid system of abelian groups $\{ A_{\alpha} \}_{\alpha \in 2^\lambda}$
such that $|A_{\alpha}|=\lambda$.
\end{thm}
\begin{prop}\label{her}
Consider a rigid system $\{ A_{\alpha} \}_{\alpha \in Y}$ of abelian groups.
If $\alpha\not=\beta$ then in $D(\mathbb Z)$,
$P_{A_{\alpha}}A_{\beta}=A_{\beta}\not=0$ hence $A_{\alpha} \not< A_{\beta}$.
\end{prop}
\proof
Suppose $\alpha\not=\beta$. Since $\{ A_{\alpha} \}$ is a rigid system,
$\Hom(A_{\alpha},A_{\beta})=0$ and since $H_i(A_{\alpha})=0$ for $i<0$ and
$H_i(A_{\beta})=0$ for $i>0$, $\Hom(s^iA_{\alpha}, A_{\beta})=0$ for all $i>0$.
Thus $P_{A_{\alpha}}A_{\beta}=A_{\beta}\not=0$.
That $A_{\alpha} \not< A_{\beta}$ then
follows directly from the definition of $<$.
\cqfd
\begin{cor}\label{properclass}
The class of t-structures, and hence also the class of
nullity classes, in $D(\mathbb Z)$ do not form a set.
\end{cor}
\proof
For any cardinal $\lambda$ let $\{ A_{\alpha} \}_{\alpha \in 2^{\lambda}}$ be
the rigid system of abelian groups of Theorem \ref{She}. To each
$A_{\alpha}$ using Theorem \ref{AJS} we associate the aisle $\overline{C(A_{\alpha})}$.
By Proposition \ref{her}, if $\alpha\not=\beta$ then
$A_{\alpha} \not< A_{\beta}$ and hence by Proposition \ref{trans}
$\overline{C(A_{\alpha})}\not=\overline{C(A_{\beta})}$. Thus the aisles
$\{ \overline{C(A_{\alpha})}\}_{\alpha\in 2^{\lambda}}$ are all distinct, which means
that the t-structures
$(\overline{C(A_{\alpha})},s\overline{C(A_{\alpha})}^{\perp})$ are also distinct. So
we see that there are at least $2^{\lambda}$
distinct t-structures. Since $\lambda$ is arbitrary the proof is complete.
\cqfd
In the category of spectra by results of Bousfield,  homological localizations are of the
form $P_A$. It was shown by Ohkawa \cite{Ohkawa}
that this subclass of localizations do form a set. It is unknown if
all localizations of the form $P_A$ which are stable under desuspension from a set.
The next theorem shows that if we take all localizations of the form $P_A$
without assuming they are closed under suspension then we do not get a set.
\begin{thm}\label{properclassinspectra}
The class of t-structures, and hence also the class of
nullity classes, in spectra do not form a set.
Similarly the class of nullity classes in spaces do not form a set.
\end{thm}
\proof
We really only give an outline of the proof.
Those initiated to the calculus of $P_A$ can easily fill in the details.
Recall that $K(G,n)$ is the Eilenberg-Mac Lane spectrum (or space) with homotopy groups
$G$ concentrated in dimension $n$. The functors $P_E$ have been constructed in spaces and spectra,
for example see \cite{DF} and \cite{Hirschhorn}. For spectra and any $E$, $\overline{C(E)}$ is an aisle.
Then for a rigid system $\{ A_{\alpha} \}_{\alpha \in 2^{\lambda}}$ of abelian groups,
all the nullity classes (and t-structures if we are working in spectra)
$\overline{C(K(A_{\alpha},n))}$ are distinct for the same reasons as above.
Since $\lambda$ is arbitrary this means that there is not a set of them.
\cqfd

%

\bigskip

\noindent\underline{Address :}\\
\noindent
Donald Stanley\\
University of Regina
\\e-mail : {\tt stanley@math.uregina.ca}
\end{document}